\input amstex
\input amsppt.sty
\magnification=\magstep1
\baselineskip=16truept
\vsize=22.2truecm
\NoRunningHeads
\nologo
\pageno=1
\TagsOnRight
\def\Z{\Bbb Z}

\def\floor #1{\left\lfloor{#1}\right\rfloor}
\def\ceil #1{\left\lceil{#1}\right\rceil}
\def\dec #1{\left\{{#1}\right\}}
\def\Proof{\noindent{\it Proof}}
\def\jacob #1#2{\bigg(\frac{#1}{#2}\bigg)}
\def\sign #1{{\roman {sign}}({#1})}
\def\Remark{\medskip\noindent{\it  Remark}}
\def\Ack{\medskip\noindent {\bf Acknowledgment}}
\topmatter
\title
A remark on Zoloterav's theorem
\endtitle
\author{Hao Pan}
\endauthor
\address
Department of Mathematics, Nanjing University,
Nanjing 210093, People's Republic of China
\endaddress
\email{haopan79\@yahoo.com.cn}\endemail
\abstract
Let $n\geq 3$ be an odd integer. For any integer $a$ prime to $n$, define the permutation $\gamma_{a,n}$ of
$\{1,\ldots,(n-1)/2\}$ by
$$
\gamma_{a,n}(x)=\cases
n-\dec{ax}_n\quad&\text{if }\dec{ax}_n\geq (n+1)/2,\\
\dec{ax}_n\quad&\text{if }\dec{ax}_n\leq (n-1)/2,\\
\endcases
$$
where $\dec{x}_n$ denotes the least nonnegative residue of $x$ modulo $n$. In this note, we show that the
sign of $\gamma_{a,n}$ coincides with the Jacobi symbol $\left(\frac{a}{n}\right)$ if $n\equiv 1\pmod{4}$,
and $1$ if $n\equiv 3\pmod{4}$.
\endabstract
\endtopmatter
\document
\heading
1. Introduction
\endheading
Let $n\geq 3$ be an odd integer. For arbitrary $x\in\Z$, let $\dec{x}_n$ denote the least nonnegative residue of $x$ modulo $n$.
Then for any integer $a$ with $(a,n)=1$, define the permutation $\sigma_{a,n}$ of $\{0,1,\ldots,n-1\} $ by $\sigma_{a,n}(x)=\dec{ax}_n$.
In his proof of quadratic reciprocity law, Zolotarev [Zo] showed that the sign of $\sigma_{a,n}$
(abbrev. as $\sign{\sigma_{a,n}}$) coincides with the Jacobi symbol $\left(\frac{a}{n}\right)$ if $n$ is prime. In fact, this result,
which is now known as Zolotarev's theorem, is also valid for composite values of $n$. The reader may find the proofs of Zolotarev's theorem
(for arbitrary odd integer $n\geq 3$) in [DS], [Me] and [Ri].
And for the sake of completeness, we include a sketch proof as a appendix.

Zolotarev's theorem actually implies the well-known Gauss lemma:
{\it
$$
\jacob{a}{n}=(-1)^{|S(a,n)|}
$$
where
$$
S(a,n)=\{1\leq x\leq (n-1)/2:\,\dec{ax}_n\geq (n+1)/2\}.
$$
}
Indeed, let
$$
S'(a,n)=\{(n+1)/2\leq x\leq n-1:\,\dec{ax}_n\leq (n-1)/2\},
$$
and define two permutations $\gamma_{a,n}$ and $\gamma_{a,n}'$ of
$\{0,1,\ldots,n-1\}$ by
$$
\gamma_{a,n}(x)=\cases
n-\dec{ax}_n\quad&\text{if }x\in S(a,n),\\
\dec{ax}_n\quad&\text{if }1\leq x\leq (n-1)/2\text{ and }x\not\in S(a,n),\\
x\quad&\text{otherwise}.\\
\endcases
$$
and
$$
\gamma_{a,n}'(x)=\cases
n-\dec{ax}_n\quad&\text{if }x\in S'(a,n),\\
\dec{ax}_n\quad&\text{if }(n+1)/2\leq x\leq n-1\text{ and }x\not\in S'(a,n),\\
x\quad&\text{otherwise}.\\
\endcases
$$
Clearly $\sign{\gamma_{a,n}}=\sign{\gamma_{a,n}'}$ since $x\in S(a,n)$ if
and only if $n-x\in S'(a,n)$. And
$$
(\gamma_{a,n}\gamma_{a,n}')(x)=\cases
n-\dec{ax}_n\quad&\text{if }x\in S(a,n)\cup S'(a,n),\\
\dec{ax}_n\quad&\text{otherwise}.\\
\endcases
$$
Thus by the Zolotarev theorem
$$
(-1)^{|S(a,n)|}=\sign{\gamma_{a,n}}\sign{\gamma_{a,n}'}\sign{\sigma_{a,n}}=\jacob{a}{n}.
$$
(In the usual statement of Gauss's lemma, $n$ is assumed to be a prime. The general version of Gauss's lemma for composite values of $n$
was proved by Jenkins [Je] and independently by Schering [Sc].)

Now a natural question arises. What is the sign of the permutation
$\gamma_{a,n}$? Observing that $\gamma_{ab,n}=\gamma_{a,n}\gamma_{b,n}$ and $\gamma_{a,n}=\gamma_{-a,n}$,
we know that $\sign{\gamma_{a,n}}$ is a real even character modulo $n$. In this note, we shall prove that
\proclaim{Theorem 1}
$$
\sign{\gamma_{a,n}}=\cases
\left(\frac{a}{n}\right)\quad&\text{if }n\equiv 1\pmod{4},\\
1\quad&\text{if }n\equiv 3\pmod{4}.\\
\endcases
$$
\endproclaim
The proof of Theorem 1 will be given in the next section.

\heading
2. Proof of Theorem 1
\endheading
Below we view $\gamma_{a,n}$ as a permutation of $\{1,\ldots,(n-1)/2\}$. Let $I(a,n)$ be the number of
inversions in $\gamma_{a,n}$. Let $\floor{\cdot}$ (resp. $\ceil{\cdot}$) denote the well-known floor function
(resp. ceiling function). And let $\dec{x}$ denote the fractional value of $x$, i.e., $\dec{x}=x-\floor{x}$.
\proclaim{Lemma 1} Suppose that $1\leq l<k\leq (n-1)/2$. Then
$$
\floor{\frac{ak}{n}}-\floor{\frac{al}{n}}-\floor{\frac{a(k-l)}{n}}=\cases
1&\text{if }\dec{ak}_n<\dec{al}_n,\\
0&\text{if }\dec{ak}_n>\dec{al}_n.\\
\endcases\tag 1
$$
And
$$
\floor{\frac{a(k+l)}{n}}-\floor{\frac{ak}{n}}-\floor{\frac{al}{n}}=\cases
1&\text{if }\dec{ak}_n+\dec{al}_n>n,\\
0&\text{if }\dec{ak}_n+\dec{al}_n<n.\\
\endcases\tag 2
$$
\endproclaim
\Proof. If $\dec{ak}_n<\dec{al}_n$, then
$$
\floor{\frac{ak}{n}}-\floor{\frac{al}{n}}=\bigg(\frac{ak}{n}-\dec{\frac{ak}{n}}\bigg)-\bigg(\frac{al}{n}-\dec{\frac{al}{n}}\bigg)>\frac{ak}{n}-\frac{al}{n}\geq \floor{\frac{a(k-l)}{n}},
$$
whence
$$
\floor{\frac{ak}{n}}=\floor{\frac{al}{n}}+\floor{\frac{a(k-l)}{n}}+1.
$$
And when $\dec{ak}_n>\dec{al}_n$, we
$$
\floor{\frac{ak}{n}}-\floor{\frac{al}{n}}<\frac{ak}{n}-\frac{al}{n}<\floor{\frac{a(k-l)}{n}}+1.
$$
So
$$
\floor{\frac{ak}{n}}=\floor{\frac{al}{n}}+\floor{\frac{a(k-l)}{n}}.
$$
If $\dec{ak}_n+\dec{al}_n>n$, then
$$
\dec{a(k+l)}_n=\dec{\dec{ak}_n+\dec{al}_n}_n=\dec{ak}_n+\dec{al}_n-n.
$$
It follows that
$$
\dec{\frac{ak}{n}}+\dec{\frac{al}{n}}=\dec{\frac{a(k+l)}{n}}+1.
$$
Thus
$$
\align
\floor{\frac{ak}{n}}+\floor{\frac{al}{n}}=&\bigg(\frac{ak}{n}-\dec{\frac{ak}{n}}\bigg)+\bigg(\frac{al}{n}-\dec{\frac{al}{n}}\bigg)\\
=&\frac{a(k+l)}{n}-\dec{\frac{a(k+l)}{n}}-1=\floor{\frac{a(k+l)}{n}}-1.
\endalign
$$
Similarly we have
$\dec{a(k+l)}_n=\dec{ak}_n+\dec{al}_n$
provided that $\dec{ak}_n+\dec{al}_n<n$. Hence
$$
\align
\floor{\frac{ak}{n}}+\floor{\frac{al}{n}}=&\bigg(\frac{ak}{n}-\dec{\frac{ak}{n}}\bigg)+\bigg(\frac{al}{n}-\dec{\frac{al}{n}}\bigg)\\
=&\frac{a(k+l)}{n}-\dec{\frac{a(k+l)}{n}}=\floor{\frac{a(k+l)}{n}}.
\endalign
$$
\qed

\Remark. (1) was used in a proof of Zolotarev's theorem by Meyer [Me] (or [RE]).

\proclaim{Lemma 2}
$$
I(a,n)\equiv\sum_{1\leq l<k\leq(n-1)/2}\bigg(\floor{\frac{a(k+l)}{n}}+\floor{\frac{a(k-l)}{n}}\bigg)\pmod{2}.
$$
\endproclaim
\Proof.
It is sufficient to show that
$$
\epsilon_{k,l}\equiv\floor{\frac{a(k+l)}{n}}+\floor{\frac{a(k-l)}{n}}\pmod{2}
$$
for $1\leq l<k\leq (n-1)/2$, where
$$
\epsilon_{k,l}=\cases
1\quad&\text{if }\gamma_{a,n}(l)>\gamma_{a,n}(k),\\
0\quad&\text{otherwise}.\\
\endcases
$$
Firstly we assume that $\epsilon_{k,l}=1$.
If $\dec{ak}_n<\dec{al}_n\leq(n-1)/2$, then
by (1), we have
$$
\floor{\frac{ak}{n}}=\floor{\frac{al}{n}}+\floor{\frac{a(k-l)}{n}}+1.
$$
And noting that $\dec{ak}_n+\dec{al}_n<2\cdot(n-1)/2$, it follows from (2) that
$$
\floor{\frac{a(k+l)}{n}}=\floor{\frac{ak}{n}}+\floor{\frac{al}{n}}.
$$
Similarly when $n-\dec{ak}_n<n-\dec{al}_n\leq(n-1)/2$, we have
$(n+1)/2\leq\dec{al}_n<\dec{ak}_n$. So applying Lemma 1,
$$
\floor{\frac{ak}{n}}=\floor{\frac{al}{n}}+\floor{\frac{a(k-l)}{n}}
$$
and
$$
\floor{\frac{a(k+l)}{n}}=\floor{\frac{ak}{n}}+\floor{\frac{al}{n}}+1.
$$
And suppose that $n-\dec{ak}_n<\dec{al}_n\leq(n-1)/2$. Then
$\dec{ak}_n+\dec{al}_n>n$
and
$\dec{al}_n\leq(n-1)/2<\dec{ak}_n$.
Thus
$$
\floor{\frac{a(k+l)}{n}}=\floor{\frac{ak}{n}}+\floor{\frac{al}{n}}+1
$$
and
$$
\floor{\frac{ak}{n}}=\floor{\frac{al}{n}}+\floor{\frac{a(k-l)}{n}}.
$$
Finally, we have
$$
\floor{\frac{a(k+l)}{n}}=\floor{\frac{ak}{n}}+\floor{\frac{al}{n}}
$$
and
$$
\floor{\frac{ak}{n}}=\floor{\frac{al}{n}}+\floor{\frac{a(k-l)}{n}}+1
$$
provided that $\dec{ak}_n<n-\dec{al}_n\leq(n-1)/2$. Thus we always obtain that
$$
\floor{\frac{ak}{n}}-\bigg(\floor{\frac{al}{n}}+\floor{\frac{a(k-l)}{n}}\bigg)
\equiv\floor{\frac{a(k+l)}{n}}-\bigg(\floor{\frac{ak}{n}}+\floor{\frac{al}{n}}\bigg)-1\pmod{2},
$$
i.e.,
$$
\floor{\frac{a(k+l)}{n}}+\floor{\frac{a(k-l)}{n}}\equiv 1+2\floor{\frac{ak}{n}}\equiv 1\pmod{2}.
$$

The arguments when $\epsilon_{k,l}=0$ are very similar, so we omit the details.\qed
\proclaim{Lemma 3}
$$
\jacob{a}{n}=(-1)^{T(a,n)+(a-1)(n^2-1)/8}
$$
where
$$
T(a,n)=\sum_{k=1}^{(n-1)/2}\floor{\frac{ak}{n}}.
$$
\endproclaim
\Proof. This lemma easily follows from Gauss's lemma and Lemma 3.1 in [Su].\qed
\Remark. When $a$ is odd, we have $\left(\frac{a}{n}\right)=(-1)^{T(a,n)}$, which is applied to a proof of quadratic reciprocity law in the most of textbooks on
number theory (cf. [Ro]).

{\medskip\noindent{\it Proof of Theorem 1}}. We compute
$$
\align
&\sum_{1\leq l<k\leq(n-1)/2}\bigg(\floor{\frac{a(k+l)}{n}}+\floor{\frac{a(k-l)}{n}}\bigg)\\
=&\sum_{k=2}^{(n-1)/2}\sum_{h=k+1}^{2k-1}\floor{\frac{ah}{n}}+\sum_{k=2}^{(n-1)/2}\sum_{h=1}^{k-1}\floor{\frac{ah}{n}}\\
=&\sum_{k=2}^{(n-1)/2}\sum_{h=1}^{2k}\floor{\frac{ah}{n}}-\sum_{k=2}^{(n-1)/2}\bigg(\floor{\frac{ak}{n}}+\floor{\frac{2ak}{n}}\bigg)\\
=&\sum_{k=1}^{(n-1)/2}\sum_{h=1}^{2k}\floor{\frac{ah}{n}}-\sum_{k=1}^{(n-1)/2}\bigg(\floor{\frac{ak}{n}}+\floor{\frac{2ak}{n}}\bigg).
\endalign
$$
Now
$$
\align
\sum_{k=1}^{(n-1)/2}\sum_{h=1}^{2k}\floor{\frac{ah}{n}}=&\sum_{h=1}^{n-1}\bigg(\frac{n-1}{2}-\ceil{\frac{h}{2}}+1\bigg)\floor{\frac{ah}{n}}\\
=&\frac{n+1}{2}\sum_{h=1}^{n-1}\floor{\frac{ah}{n}}-\sum_{h=1}^{n-1}\ceil{\frac{h}{2}}\floor{\frac{ah}{n}}.
\endalign
$$
Observe that
$$
\sum_{h=1}^{n-1}\floor{\frac{ah}{n}}=\sum_{h=1}^{n-1}\bigg(\frac{ah}{n}-\dec{\frac{ah}{n}}\bigg)=(a-1)\sum_{h=1}^{n-1}\frac{h}{n}=\frac{(a-1)(n-1)}{2}.
$$
And
$$
\align
\sum_{h=1}^{n-1}\ceil{\frac{h}{2}}\floor{\frac{ah}{n}}
=&\sum_{k=1}^{(n-1)/2}k\bigg(\floor{\frac{2ak}{n}}+\floor{\frac{a(2k-1)}{n}}\bigg)\\
=&\sum_{k=1}^{(n-1)/2}k\bigg(\frac{2ak}{n}+\frac{a(2k-1)}{n}-\dec{\frac{2ak}{n}}-\dec{\frac{a(2k-1)}{n}}\bigg)\\
=&\frac{a(4n-3)(n+1)(n-1)}{24n}-\sum_{k=1}^{(n-1)/2}k\bigg(\dec{\frac{2ak}{n}}+\dec{\frac{a(2k-1)}{n}}\bigg).
\endalign
$$
Finally,
$$
\align
&\sum_{k=1}^{(n-1)/2}k\dec{\frac{2ak}{n}}+\sum_{k=1}^{(n-1)/2}k\dec{\frac{a(2k-1)}{n}}\\
=&\sum_{k=1}^{(n-1)/2}k\dec{\frac{2ak}{n}}+\sum_{k=1}^{(n-1)/2}k\bigg(1-\dec{\frac{a(n-2k+1)}{n}}\bigg)\\
=&\sum_{k=1}^{(n-1)/2}k\dec{\frac{2ak}{n}}+\sum_{k=1}^{(n-1)/2}\bigg(\frac{n+1}{2}-k\bigg)\bigg(1-\dec{\frac{2ak}{n}}\bigg)\\
=&2\sum_{k=1}^{(n-1)/2}k\dec{\frac{2ak}{n}}+\frac{n+1}{2}\sum_{k=1}^{(n-1)/2}\bigg(1-\dec{\frac{2ak}{n}}\bigg)-\sum_{k=1}^{(n-1)/2}k\\
\equiv&\frac{n^2-1}{8}-\frac{n+1}{2}\sum_{k=1}^{(n-1)/2}\bigg(\frac{2ak}{n}-\floor{\frac{2ak}{n}}\bigg)\\
=&\frac{n^2-1}{8}-\frac{a(n+1)^2(n-1)}{8n}+\frac{n+1}{2}\sum_{k=1}^{(n-1)/2}\floor{\frac{2ak}{n}}\pmod{2}.
\endalign
$$
Since $n$ is odd,
$$
\frac{n+1}{2}\cdot\frac{(a-1)(n-1)}{2}\equiv\frac{a(n+1)^2(n-1)}{8n}\equiv 0\pmod{2}
$$
and
$$
\frac{a(4n-3)(n+1)(n-1)}{24n}\equiv\frac{a(n^2-1)}{8}\pmod{2}.
$$
Thus applying Lemma 2, we have
$$
I(a,n)\equiv\frac{(1-a)(n^2-1)}{8}+\frac{n-1}{2}\sum_{k=1}^{(n-1)/2}\floor{\frac{2ak}{n}}-\sum_{k=1}^{(n-1)/2}\floor{\frac{ak}{n}}\pmod{2}.
$$
Therefore by Lemma 3, when $n\equiv 1\pmod{4}$,
$$
\sign{\gamma_{a,n}}=(-1)^{I(a,n)}=(-1)^{\sum_{k=1}^{(n-1)/2}\floor{ak/n}+(a-1)(n^2-1)/8}=\jacob{a}{n}.
$$
And if $n\equiv 3\pmod{4}$, then
$$
\align
\sign{\gamma_{a,n}}=(-1)^{I(a,n)}=&(-1)^{\sum_{k=1}^{(n-1)/2}(\floor{ak/n}+\floor{2ak/n})+(a-1)(n^2-1)/8}\\
=&(-1)^{(2a-1)(n^2-1)/8}\jacob{a}{n}\jacob{2a}{n}=1
\endalign
$$
since $(-1)^{(n^2-1)/8}=\left(\frac{2}{n}\right)$. We are done.\qed

\heading
3. Appendix: a sketch proof of Zolotarev's theorem
\endheading
Clearly $\sigma_{ab,n}=\sigma_{a,n}\sigma_{b,n}$.
(i) Assume that $n=p$ where $p$ is a prime. Let $g$ be a prime root of $p$. Then it suffices to show that
$\sign{\sigma_{g,p}}=\left(\frac{g}{p}\right)=-1$. Clearly the action of $\sigma_{g,p}$ on $\{1,2,\ldots,p-1\}$ coincides with the $(p-1)$-cycle $(\{g^{p-2}\}_p\ \{g^{p-3}\}_p\ \ldots\ \{g\}_p\ 1)$.
Hence $\sigma_{g,p}$ is an odd permutation. (ii) Suppose that $n=p^\alpha$ where $p$ is a prime and $\alpha\geq 2$. We need the following simple fact:

\medskip{\it Let $\sigma$ be a permutation of the finite set $X$. And suppose that $Y$ be a nonempty proper subset of $X$ such that $Y$ is $\sigma$-invariant, i.e.,
$\sigma(Y)=Y$. Then $\sign{\sigma}=\sign{\sigma|_{Y}}\sign{\sigma|_{X\setminus Y}}$,
where $\sigma|_{Y}$ denotes the restriction of $\sigma$ to $Y$.}

{\medskip\noindent}Now let $X_1=\{0,p,2p,\ldots,p(p^{\alpha-1}-1)\}$ and $X_2=\{0,1,\ldots,p^{\alpha}-1\}\setminus X_1$. Let $g$ be a prime root of
$p^{\alpha}$. Then following the argument in (i), we have $\sign{\sigma_{g,p^{\alpha}}|_{X_2}}=-1$. And note that $\sigma_{g,p^{\alpha}}|_{X_1}$
is equivalent to $\sigma_{g,p^{\alpha-1}}$. Thus by an induction on $\alpha$,
$$
\sign{\sigma_g}=\sign{\sigma_g|_{X_1}}\sign{\sigma_g|_{X_2}}=\jacob{g}{p^{\alpha-1}}\cdot(-1)=\jacob{g}{p^\alpha}.
$$
(iii) If $n=p_1^{\alpha_1}p_2^{\alpha_2}\cdots p_k^{\alpha_k}$ where these $p_i$ are distinct primes and $\alpha_i\geq 1$,
then by the Chinese remainder theorem we can choose $g_i$ for $1\leq i\leq k$ such that $g_i$ is a prime root
of $p_i^{\alpha_i}$ and $g_i\equiv 1\pmod{n/p_i^{\alpha_i}}$. Let
$$
X_{a_1,a_2,\ldots,a_{k-1}}=\{0\leq x<n:\, x\equiv a_i\pmod{p^{\alpha_i}}\text{ for }1\leq i\leq k-1\}.
$$
It is easy to see that $\sigma_{g_k,n}(X_{a_1,\ldots,a_{k-1}})=X_{a_1,\ldots,a_{k-1}}$ and
$\bigcup_{\Sb 0\leq a_i<p_{i}^{\alpha_i}\\1\leq i\leq k-1\endSb}X_{a_1,\ldots,a_{k-1}}=\{0,1,\ldots,n-1\}$.
Then by (ii), we have
$$
\sign{\sigma_{g_k,n}|_{X_{a_1,\ldots,a_{k-1}}}}=\sign{\sigma_{g_k,p_k^{\alpha_k}}}=\jacob{g_k}{p_k^{\alpha_k}}=\jacob{g_k}{n}
$$
for arbitrary $a_1,\ldots,a_{k-1}$. Therefore
$$
\sign{\sigma_{g_k,n}}=\prod_{\Sb 0\leq a_i<p_{i}^{\alpha_i}\\1\leq i\leq k-1\endSb}\sign{\sigma_{g_k,n}|_{X_{a_1,\ldots,a_{k-1}}}}
=\jacob{g_k}{n}^{n/p_k^{\alpha_k}}=\jacob{g_k}{n}.
$$
Similarly we have $\sign{\sigma_{g_i,n}}=\left(\frac{g_i}{n}\right)$ for each $1\leq i\leq k$. This concludes our proof since
$g_1, g_2,\ldots, g_k$ generate the modulo multiplication group of $n$.
\widestnumber\key{AA}

\Ack. I thank my advisor, Zhi-Wei Sun, for his help on this paper.

\enddocument